\documentclass[12pt]{article}%
\usepackage{amsmath}
\usepackage{amsfonts}
\usepackage[twocolumn=false]{geometry}
\usepackage{graphicx}%
\usepackage{amssymb}
\providecommand{\U}[1]{\protect\rule{.1in}{.1in}}

\begin{document}

\title{On a Subclass of p-Valent Functions with Negative \ Coefficients Defined by
Using Rafid Operator}
\author{A. H. El-Qadeem$^{\text{1}}$\ \ \ \& \ S. K. Al-ghazal$^{\text{2}}$
\and Department of Mathematics, Faculty of Science, Zagazig University,\\Zagazig 44519, Al-Sharkia, Egypt
\and ahhassan@science.zu.edu.eg$^{\text{1}}$ \ \& \ z1z10z78z@gmail.com$^{\text{2}%
}$}
\date{}
\maketitle

\begin{abstract}
By using Rafid operator we define the subclass $R_{\mu,p}^{\delta}%
(\alpha;A,B)\ $and \ \ \ $P_{\mu,p}^{\delta}(\alpha;A,B)$ of analytic and
p-valent functions with negative coefficients we investigate some sharp
results including coefficients estimates, distortion theorem, radii of
starlikeness, convexity, close-to-convexity, and modified-Hadamard product.
Finally, we give an application of fractional calculus and
Bernadi-Libora-Livingstion operator.

\end{abstract}



\section*{1. Introduction}

Let $T(p)$ denotes the class of normalized p-valent functions $f$ which are
analytic in $U$ =$\{z\in%
\mathbb{C}
:|z|$ $<1\}$, and given\ by%
\begin{equation}
f(z)=z^{p}-%
{\textstyle\sum\limits_{k=p+1}^{\infty}}
a_{k}z^{k}\text{ \ }(a_{k}\geq0,p\in%
\mathbb{N}
=\{1,2,3...\}),\tag{1.1}%
\end{equation}
A function $f\in T(p)$ is called p-valent starlike of order $\alpha
(0\leq\alpha<p),$ if and only if%
\begin{equation}
\operatorname*{Re}\left\{  \frac{zf^{^{\prime}}(z)}{f(z)}\right\}
>\alpha\text{ \ \ }\left(  0\leq\alpha<p;z\in U\right)  ,\tag{1.2}%
\end{equation}
\noindent we denote by $T^{\ast}(p,\alpha)$ the class of all p-valent starlike
functions of order $\alpha$. Also a function $f\in T(p)$ is called p-valent
convex of order $\alpha(0\leq\alpha<p),$ if and only if%
\begin{equation}
\operatorname*{Re}\left\{  1+\frac{zf^{^{\prime\prime}}(z)}{f^{^{\prime}}%
(z)}\right\}  >\alpha\text{ \ \ }\left(  0\leq\alpha<p;z\in U\right)
,\tag{1.3}%
\end{equation}
we denote by $C(p,\alpha)$ the class of all p-valenty convex functions of
order $\alpha$. For more informations about the subclasses $T^{\ast}%
(p,\alpha)$ and $C(p,\alpha)$, see \cite{Sliverman}.\newline Motivated by
Atshan and Rafid see \cite{sp 1}, we introduce the following p-valent analogue
$R_{\mu,p}^{\delta}:T(p)\longrightarrow T(p):$\newline For $0\leq\mu<1$ and
$0\leq\delta\leq1,$%

\begin{equation}
R_{\mu,p}^{\delta}f(z)=\frac{1}{\Gamma(p+\delta)(1-\mu)^{p+\delta}}%
{\displaystyle\int\limits_{0}^{\infty}}
t^{\delta-1}e^{-\left(  \frac{t}{1-\mu}\right)  }f(zt)dt\tag{1.4}%
\end{equation}

\noindent where $\Gamma$\ stands for Euler's Gamma function (which is valid
for all complex numbers except the non-positive integers).

\noindent Let $f\ $and\ $g$\ be analytic in $U$. Then we say that the function
$g$\ is subordinate to$\ f$\ if there exists an analytic function\ in $U$ such
that\ $\left\vert w(z)\right\vert <1$ $(z\in U)$\ and\ $g(z)=f(w(z))$. For
this subordination, the symbol\ $g(z)\prec f(z)$ is used. In case $f(z)$ is
univalent in $U$, the subordination $g(z)\prec f(z)$ is equivalent
to\ $g(0)=f(0)$ and $\ g(U)\subset f(U),$ see \cite{Miller and Mocanu}.

\noindent For the function $f$ given by (1.1) and $g(z)=z^{p}-%
{\textstyle\sum\nolimits_{k=p+1}^{\infty}}
b_{k}z^{k}$, the modified Hadamard product (or convolution) of $f$ and $g$ is
denoted by $f\ast g$ and is given by%
\begin{equation}
\left(  f\ast g\right)  (z)=z^{p}-%
{\textstyle\sum\limits_{k=p+1}^{\infty}}
a_{k}b_{k}z^{k}=(g\ast f)(z).\tag{1.5}%
\end{equation}

\noindent\textbf{Definition 1. }For $-1\leq B<A\leq1$ and $0\leq\alpha<p$, let
$R_{\mu,p}^{\delta}(\alpha;A,B)$\ be the subclass of functions $f\in T(p)$ for which:%

\begin{equation}
\frac{z(R_{\mu,p}^{\delta}f(z))^{^{\prime}}}{R_{\mu,p}^{\delta}f(z)}%
\prec(p-\alpha)\frac{1+Az}{1+Bz}+\alpha,\tag{1.6}%
\end{equation}
that is, that
\begin{equation}
R_{\mu,p}^{\delta}(p,\alpha;A,B)=\left\{  f\in T(p):\left\vert \tfrac
{\tfrac{z(R_{\mu,p}^{\delta}f(z))^{\prime}}{R_{\mu,p}^{\delta}f(z)}-p}%
{B\tfrac{z(R_{\mu,P}^{\delta}f(z))^{\prime}}{R_{\mu,p}^{\delta}f(z)}-\left[
Bp+(A-B)(p-\alpha)\right]  }\right\vert <1,z\in U\right\}  .\tag{1.7}%
\end{equation}
Note that $\operatorname{Re}\left\{  (p-\alpha)\frac{1+Az}{1+Bz}%
+\alpha\right\}  >\tfrac{1-A+\alpha(A-B)}{1-B}.$\newline Also, for $-1\leq
B<A\leq1$ and $0\leq\alpha<p$, let $P_{\mu,p}^{\delta}(\alpha;A,B)$ be
subclass of functions $f$ $\in T$ $(p)$ for which:%

\begin{equation}
1+\frac{z(R_{\mu,p}^{\delta}f(z))^{^{\prime\prime}}}{(R_{\mu,p}^{\delta
}f(z))^{^{\prime}}}\prec(p-\alpha)\frac{1+Az}{1+Bz}+\alpha,\tag{1.8}%
\end{equation}
For (1.6) and (1.8) it is clear that%

\begin{equation}
f(z)\in P_{\mu,p}^{\delta}\iff\frac{zf^{^{\prime}}(z)}{p}\in R_{\mu,p}%
^{\delta}\tag{1.9}%
\end{equation}

\section*{2. Main Results}

Unless otherwise mentioned, we assume in the reminder of this paper that,
$0\leq\alpha<p,$ $0\leq$ $\mu<1$, $0\leq\delta\leq1,$

$-1\leq B<A\leq1,$ $p\in%
\mathbb{N}
$ and $z\in U.$

\subsection*{2.1. Coefficients Estimate}

\noindent\textbf{Theorem 1}. Let the function $f(z)$ be given by (1.1). Then
$f(z)\in R_{\mu,p}^{\delta}(\alpha;A,B),$ if and only if%

\begin{equation}%
{\textstyle\sum\limits_{k=p+1}^{\infty}}
\left[  (1-B)(k-p)+(A-B)(p-\alpha)\right]  (1-\mu)^{k-p}\frac{\Gamma
(k+\delta)}{\Gamma(p+\delta)}a_{k}\leq(A-B)(p-\alpha).\tag{2.1}%
\end{equation}
\noindent\textbf{Proof}. Assume that the inequality (2.1) holds true. We find
from (1.1) and (2.1) that

\noindent thus we have%

\[
\left\vert z(R_{\mu,p}^{\delta}f(z))^{^{\prime}}-p(R_{\mu,p}^{\delta
}f(z))\right\vert -\left\vert B\left[  z(R_{\mu,p}^{\delta}f(z))^{^{\prime}%
}\right]  -[Bp+(A-B)(p-\alpha)][R_{\mu,p}^{\delta}f(z)]\right\vert
\]

\[%
{\textstyle\sum\limits_{k=p+1}^{\infty}}
[(1-B)(k-p)+(A-B)(p-\alpha)]\frac{\Gamma(k+\delta)}{\Gamma(p+\delta)}%
(1-\mu)^{k-p}a_{k}-(A-B)(p-\alpha)\leq0.
\]
\noindent Hence, by the maximum modulus theorem, we have%

\[
\left\vert \frac{\dfrac{z(R_{\mu,p}^{\delta}f(z))^{^{\prime}}}{R_{\mu
,p}^{\delta}f(z)}-p}{B\dfrac{z(R_{\mu,p}^{\delta}f(z))^{^{\prime}}}{R_{\mu
,p}^{\delta}f(z)}-\left[  Bp+(A-B)(p-\alpha)\right]  }\right\vert <1.
\]

\noindent Thus $f\in R_{\mu,p}^{\delta}(\alpha;A,B).$

\noindent Conversely, let $f\in R_{\mu,p}^{\delta}(\alpha;A,B)$ be given by
(1.1), then from (1.1) and (1.7), we have%

\[
\left\vert \frac{\dfrac{z(R_{\mu,p}^{\delta}f(z))^{^{\prime}}}{R_{\mu
,p}^{\delta}f(z)}-p}{B\dfrac{z(R_{\mu,p}^{\delta}f(z))^{^{\prime}}}{R_{\mu
,p}^{\delta}f(z)}-\left[  Bp+(A-B)(p-\alpha)\right]  }\right\vert
\]

\[
=\left\vert \frac{%
{\textstyle\sum\limits_{k=p+1}^{\infty}}
(k-p)(1-\mu)^{k-p}\dfrac{\Gamma(k+\delta)}{\Gamma(p+\delta)}a_{k}z^{k}%
}{(A-B)(p-\alpha)z^{p}+%
{\textstyle\sum\limits_{k=p+1}^{\infty}}
[-B(k-p)+(A+B)(p-\alpha)](1-\mu)^{k-p}\dfrac{\Gamma(k+\delta)}{\Gamma
(p+\delta)}a_{k}z^{k}}\right\vert <1.
\]
Since $\operatorname*{Re}(z)\leq\left\vert z\right\vert $ for all $z$, we have%
\begin{equation}
\operatorname*{Re}\left\{  \frac{%
{\textstyle\sum\limits_{k=p+1}^{\infty}}
(k-p)(1-\mu)^{k-p}\dfrac{\Gamma(k+\delta)}{\Gamma(p+\delta)}a_{k}z^{k}%
}{(A-B)(p-\alpha)z^{p}+%
{\textstyle\sum\limits_{k=p+1}^{\infty}}
[-B(k-p)+(A+B)(p-\alpha)](1-\mu)^{k-p}\dfrac{\Gamma(k+\delta)}{\Gamma
(p+\delta)}a_{k}z^{k}}\right\}  <1,\tag{2.2}%
\end{equation}

\noindent choose values of $z$ on the real axis so that $\dfrac{z(R_{\mu
,p}^{\delta}f(z))^{^{\prime}}}{R_{\mu,p}^{\delta}f(z)}$ is real. It is
clearing the denominator in (2.2) and letting $z\rightarrow1^{-}$ through real
values, we have%

\begin{equation}%
{\textstyle\sum\limits_{k=p+1}^{\infty}}
(k-p)(1-\mu)^{k-p}\tfrac{\Gamma(k+\delta)}{\Gamma(p+\delta)}a_{k}%
\leq(A-B)(p-\alpha)-%
{\textstyle\sum\limits_{k=p+1}^{\infty}}
[-B(k-p)+(A-B)(p-\alpha)](1-\mu)^{k-p}\tfrac{\Gamma(k+\delta)}{\Gamma
(p+\delta)}a_{k}.\tag{2.3}%
\end{equation}

\noindent This gives the required condition.

\noindent\textbf{Corollary 1}. Let the function $f(z)$ defined by (1.1) be in
the class $R_{\mu,p}^{\delta}(\alpha;A,B)$. Then%

\begin{equation}
a_{k}\leq\frac{(A-B)(p-\alpha)}{[(1-B)(k-p)+(A-B)(p-\alpha)](1-\mu)^{k-p}%
\frac{\Gamma(k+\delta)}{\Gamma(p+\delta)}}\text{ \ \ \ \ \ }(k\geq
p+1),\tag{2.4}%
\end{equation}

\noindent the result is sharp for the function $f_{0}$ given by%

\begin{equation}
f_{0}(z)=z^{p}-\frac{(A-B)(p-\alpha)}{[(1-B)(k-p)+(A-B)(p-\alpha
)](1-\mu)^{k-p}\frac{\Gamma(k+\delta)}{\Gamma(p+\delta)}}z^{k}\text{
\ \ }(k\geq p+1).\tag{2.5}%
\end{equation}
By using (1.9) and Theorem 1, it is easily to obtain the following
result.\newline\textbf{Theorem 2}. Let the function $f(z)$ be given by (1.1).
Then $f\in P_{\mu,p}^{\delta}(\alpha;A,B)$ if and only if%

\begin{equation}%
{\textstyle\sum\limits_{k=p+1}^{\infty}}
\frac{k}{p}[(k-p)(1-B)+(A-B)(p-\alpha)](1-\mu)^{k-p}\frac{\Gamma(k+\delta
)}{\Gamma(p+\delta)}a_{k}\leq(A-B)(1-\alpha)\frac{\Gamma(k+\delta)}%
{\Gamma(p+\delta)}\tag{2.6}%
\end{equation}
\bigskip

\noindent\textbf{Corollary 2.} Let the function $f(z)$\ defined by (1.1) be in
the class $P_{\mu,p}^{\delta}(\alpha;A,B)$. Then%

\begin{equation}
a_{k}\leq\frac{(A-B)(p-\alpha)}{\frac{k}{p}[(1-B)(k-p)+(A-B)(p-\alpha
)](1-\mu)^{k-p}\frac{\Gamma(k+\delta)}{\Gamma(p+\delta)}}\text{ \ \ \ }(k\geq
p+1),\tag{2.7}%
\end{equation}

\noindent the result is sharp for the function $f_{1}(z)$ give by%

\[
f_{1}(z)=z^{p}-\frac{(A-B)(p-\alpha)}{\frac{k}{p}[(1-B)(k-p)+(A-B)(p-\alpha
)](1-\mu)^{k-p}\frac{\Gamma(k+\delta)}{\Gamma(p+\delta)}}z^{k}\text{
\ \ \ }(k\geq p+1),
\]

\subsection*{2.2. Distortion Theorem}

\noindent\textbf{Theorem 3}. If the function $f(z)$ defined by (1.1) is in the
class $R_{\mu,p}^{\delta}(\alpha;A,B).$ Then%
\[
\left[  \delta(p,m)-\frac{(A-B)(p-\alpha)(p+1)!}{(p+1-m)![(1-B)+(A-B)(p-\alpha
)](1-\mu)\frac{\Gamma(p+\delta+1)}{\Gamma(p+\delta)}}\left\vert z\right\vert
\right]  \left\vert z\right\vert ^{p-m}%
\]%
\[
\leq\left\vert f^{(m)}(z)\right\vert \leq
\]%
\begin{equation}
\left[  \delta(p,m)+\frac{(A-B)(p-\alpha)(p+1)!}{(p+1-m)![(1-B)+(A-B)(p-\alpha
)](1-\mu)\frac{\Gamma(p+\delta+1)}{\Gamma(p+\delta)}}\left\vert z\right\vert
\right]  \left\vert z\right\vert ^{p-m}.\tag{2.6}%
\end{equation}
$(m\in%
\mathbb{N}
_{0;\text{ }}p>\{m\})$

\noindent The result is sharp for the function $f_{0}$ given by%
\[
f_{0}(z)=z^{p}-\frac{(A-B)(p-\alpha)}{[(1-B)+(A-B)(p-\alpha)](1-\mu
)\frac{\Gamma(p+\delta+1)}{\Gamma(p+\delta)}}z^{p+1}%
\]

\noindent\textbf{Proof}. In view of Theorem 1, we have%
\[
\frac{\lbrack(1-B)+(A-B)(p-\alpha)](1-\mu)\frac{\Gamma(p+\delta+1)}%
{\Gamma(p+\delta)}}{(A-B)(p-\alpha)(p+1)!}%
{\textstyle\sum\limits_{k=p+1}^{\infty}}
k!a_{k}\leq
\]%
\[%
{\textstyle\sum\limits_{k=p+1}^{\infty}}
\frac{[(1-B)(k-p)+(A-B)(p-\alpha)](1-\mu)^{k-p}\frac{\Gamma(k+\delta)}%
{\Gamma(p+\delta)}}{(A-B)(p-\alpha)(p+1)!}a_{k}\leq1,
\]
\noindent which readily yields%
\begin{equation}%
{\textstyle\sum\limits_{k=p+1}^{\infty}}
k!a_{k}\leq\frac{(A-B)(p-\alpha)(p+1)!}{[(1-B)+(A-B)(p-\alpha)](1-\mu
)\frac{\Gamma(p+\delta+1)}{\Gamma(p+\delta)}}.\tag{2.7}%
\end{equation}
\noindent Now by differentiating both sides of (1.1) m-times we have%

\begin{equation}
f^{(m)}(z)=\delta(p,m)z^{p-m}-%
{\textstyle\sum\limits_{k=p+1}^{\infty}}
\delta(k,m)a_{k}z^{k-m}.\tag{2.8}%
\end{equation}
and Theorem 3 would follow from (2.7) and (2.8).

\subsection*{2.3. Radii of Starlikeness, Convexity and Close-to-Convexity}

\noindent\textbf{Theorem 4}. Let the function $f(z)$ defined by (1.1) be in
the class $R_{\mu,p}^{\delta}(\alpha;A,B),$ then

(i) $f(z)$ is p-valently starlike of order $\zeta(0\leq\zeta<p)$ in
$\left\vert z\right\vert <r_{1},$ where%

\begin{equation}
r_{1}=\underset{k}{\inf}\left[  \frac{[(1-B)(k-p)+(A-B)(p-\alpha
)](1-\mu)^{k-p}\frac{\Gamma(k+\delta)}{\Gamma(p+\delta)}}{(A-B)(p-\alpha
)}\left(  \frac{p-\zeta}{k-\zeta}\right)  \right]  ^{\frac{1}{k-p}}(k\geq
p+1),\tag{2.9}%
\end{equation}

(ii) $f(z)$ is p-valently convex of order $\zeta(0\leq\zeta<p)$ in $\left\vert
z\right\vert <r_{2},$ where%

\begin{equation}
r_{2}=\underset{k}{\inf}\left[  \frac{[(1-B)(k-p)+(A-B)(p-\alpha
)](1-\mu)^{k-p}\frac{\Gamma(k+\delta)}{\Gamma(p+\delta)}}{(A-B)(p-\alpha
)}\left(  \frac{p(p-\zeta}{k(k-\zeta)}\right)  \right]  ^{\frac{1}{k-p}}(k\geq
p+1),\tag{2.10}%
\end{equation}
(iii) $f(z)$ is p-valently close-to-convex of order $\zeta(0\leq\zeta<p)$ in
$\left\vert z\right\vert <r_{3\text{ }}$where%
\begin{equation}
r_{3}=\underset{k}{\inf}\left[  \frac{[(1-B)(k-p)+(A-B)(p-\alpha
)](1-\mu)^{k-p}\frac{\Gamma(k+\delta)}{\Gamma(p+\delta)}}{(A-B)(p-\alpha
)}\left(  \frac{(p-\zeta}{k}\right)  \right]  ^{\frac{1}{k-p}}(k\geq
p+1),\tag{2.11}%
\end{equation}
Each of these results are sharp for the function $f(z)$ given by (2.5)

\noindent\textbf{Proof}. It is sufficient to show that%

\begin{equation}
\left\vert \frac{zf^{^{\prime}}(z)}{f(z)}-p\right\vert \leq p-\zeta\text{
\ \ \ }(\left\vert z\right\vert <r_{1};0\leq\zeta<p),\tag{2.12}%
\end{equation}
or%

\begin{equation}
\left\vert \frac{zf^{^{\prime}}(z)}{f(z)}-p\right\vert =\left\vert \frac{%
{\textstyle\sum\limits_{k=p+1}^{\infty}}
(k-p)a_{k}z^{k-p}}{z^{p}-%
{\textstyle\sum\limits_{k=p+1}^{\infty}}
a_{k}z^{k-p}}\right\vert \leq\frac{%
{\textstyle\sum\limits_{k=p+1}^{\infty}}
(k-p)a_{k}\left\vert z\right\vert ^{k-p}}{1-%
{\textstyle\sum\limits_{k=p+1}^{\infty}}
a_{k}\left\vert z\right\vert ^{k-p}}.\tag{2.13}%
\end{equation}

\noindent Inequality (2.12) holds true, when%

\[
\frac{%
{\textstyle\sum\limits_{k=p+1}^{\infty}}
(k-p)a_{k}\left\vert z\right\vert ^{k-p}}{1-%
{\textstyle\sum\limits_{k=p+1}^{\infty}}
a_{k}\left\vert z\right\vert ^{k-p}}\leq p-\zeta,
\]

\noindent or, when%

\begin{equation}%
{\textstyle\sum\limits_{k=p+1}^{\infty}}
\left(  \frac{k-\zeta}{p-\zeta}\right)  a_{k}\left\vert z\right\vert
^{k-p}\leq1,\tag{2.14}%
\end{equation}

\noindent using inequality (2.1), then (2.14) holds true if%

\begin{equation}
\left(  \frac{k-\zeta}{p-\zeta}\right)  a_{k}\left\vert z\right\vert
^{k-p}\leq\frac{\lbrack(1-B)(k-p)+(A-B)(p-\alpha)](1-\mu)^{k-p}\frac
{\Gamma(k+\delta)}{\Gamma(p+\delta)}}{(A-B)(p-\alpha)}a_{k},\text{ \ \ }(k\geq
p+1),\tag{2.15}%
\end{equation}

\noindent or%

\begin{equation}
\left\vert z\right\vert \leq\left\{  \frac{\lbrack(1-B)(k-p)+(A-B)(p-\alpha
)](1-\mu)^{k-p}\frac{\Gamma(k+\delta)}{\Gamma(p+\delta)}}{(A-B)(p-\alpha
)}\left(  \frac{p-\zeta}{k-\zeta}\right)  \right\}  ^{\frac{1}{k-p}}%
\text{\ }\ \ \text{\ }(k\geq p+1),\tag{2.16}%
\end{equation}

\noindent or%

\begin{equation}
r_{1}=\underset{k}{\inf}\left\{  \frac{[(1-B)(k-p)+(A-B)(p-\alpha
)](1-\mu)^{k-p}\frac{\Gamma(k+\delta)}{\Gamma(p+\delta)}}{(A-B)(p-\alpha
)}\left(  \frac{p-\zeta}{k-\zeta}\right)  \right\}  ^{\frac{1}{k-p}}%
\text{\ }\ \ \text{\ }(k\geq p+1).\tag{2.17}%
\end{equation}

\noindent This completes the proof (2.9).

\noindent To prove (ii) and (iii) it is sufficient to note that%

\begin{equation}
\left\vert 1+\frac{zf^{^{^{\prime\prime}}}(z)}{f^{^{\prime}}(z)}-p\right\vert
\leq p-\zeta\text{ \ \ \ }(\left\vert z\right\vert <r_{2};0\leq\zeta
<p)\tag{2.18}%
\end{equation}
and%

\begin{equation}
\left\vert \frac{f^{^{^{\prime}}}(z)}{z^{p-1}}-p\right\vert \leq p-\zeta\text{
\ \ \ }(\left\vert z\right\vert <r_{3};0\leq\zeta<p)\tag{2.19}%
\end{equation}
respectively.

\subsection*{2.4 Modified-Hadamard Product}

\noindent In this subsection, we obtain some results of the modified Hadamard
product of functions $f_{1}$ and $f_{2},$ which are defined by%

\begin{equation}
f_{v}(z)=z^{p}-%
{\textstyle\sum\limits_{k=p+1}^{\infty}}
a_{k,\nu}z^{k}\text{ \ \ \ }(a_{k,v}\geq0,v=1,2),\tag{2.20}%
\end{equation}

\noindent\textbf{Theorem 5. }If $f_{v}\in R_{\mu,p}^{\delta}(\alpha;A,B)$
$(v=1,2)$ defined by (2.20), then $(f_{1}\ast f_{2})(z)\in R_{\mu,p}^{\delta
}(\lambda;A,B)$, where%

\begin{equation}
\lambda=p-\frac{(1-B)(A-B)(p-\alpha)^{2}}{[(1-B)+(A-B)(p-\alpha)]^{2}%
(1-\mu)\tfrac{\Gamma(p+\delta+1)}{\Gamma(p+\delta)}-[(A-B)(p-\alpha)]^{2}%
},\tag{2.21}%
\end{equation}

\noindent The result is sharp for that function $f_{v}(z)$ $(v=1,2)$ given by%

\begin{equation}
f_{v}(z)=z^{p}-\frac{(A-B)(p-\alpha)}{[(1-B)+(A-B)(p-\alpha)](1-\mu
)\tfrac{\Gamma(p+\delta+1)}{\Gamma(p+\delta)}}z^{p+1}.\tag{2.22}%
\end{equation}

\noindent\textbf{Proof}. Employing the technique used earlier by Schild and
Silverman \cite{Schild and Silverman}, we need to find the largest $\lambda$
such that%

\begin{equation}%
{\textstyle\sum\limits_{k=p+1}^{\infty}}
\frac{[(1-B)(k-p)+(A-B)(p-\lambda)](1-\mu)^{k-p}\dfrac{\Gamma(k+\delta
)}{\Gamma(p+\delta)}}{(A-B)(p-\lambda)}a_{k,1}a_{k,2}\leq1\text{.}\tag{2.23}%
\end{equation}

\noindent Since $f_{v}\in R_{\mu,p}^{\delta}(p,\alpha;A,B)$ $(v=1,2)$, we
readily see that%
\[%
{\textstyle\sum\limits_{k=p+1}^{\infty}}
\frac{[(1-B)(k-p)+(A-B)(p-\alpha)](1-\mu)^{k-p}\dfrac{\Gamma(k+\delta)}%
{\Gamma(p+\delta)}}{(A-B)(p-\alpha)}a_{k,v}\leq1\ (v=1,2).
\]

\noindent Therefore, by the Cauchy-Schwarz inequality, we obtain%

\begin{equation}%
{\textstyle\sum\limits_{k=p+1}^{\infty}}
\frac{[(1-B)(k-p)+(A-B)(p-\alpha)](1-\mu)^{k-p}\dfrac{\Gamma(k+\delta)}%
{\Gamma(p+\delta)}}{(A-B)(p-\alpha)}\sqrt{a_{k,1}a_{k,2}}\leq1.\text{
}\tag{2.24}%
\end{equation}

\noindent From (2.23) and (2.24), we need only to show that%
\[
\tfrac{\lbrack(1-B)(k-p)+(A-B)(p-\lambda)]}{(p-\lambda)}\sqrt{a_{k,1}a_{k,2}%
}\leq\tfrac{\lbrack(1-B)(k-p)+(A-B)(p-\alpha)]}{(p-\alpha)}a_{k,1}%
a_{k,2}\text{ \ }(k\geq p+1),
\]
\noindent or, equivalently, that%
\begin{equation}
\sqrt{a_{k,1}a_{k,2}}\leq\frac{(p-\lambda)[(1-B)(k-p)+(A-B)(p-\alpha
)]}{(p-\alpha)[(1-B)(k-p)+(A-B)(p-\lambda)]}\text{ \ }(k\geq p+1).\tag{2.25}%
\end{equation}
\noindent Hence, in the light of inequality (2.24). It is sufficient to prove
that%
\[
\frac{(A-B)(p-\alpha)}{[(1-B)(k-p)+(A-B)(p-\alpha)](1-\mu)^{k-p}\dfrac
{\Gamma(k+\delta)}{\Gamma(p+\delta)}}%
\]%
\begin{equation}
\leq\frac{(p-\lambda)[(1-B)(k-p)+(A-B)(p-\alpha)]}{(p-\alpha
)[(1-B)(k-p)+(A-B)(p-\lambda)]}\text{ \ \ }(k\geq p+1).\tag{2.26}%
\end{equation}
\noindent It follows from (2.26) that,%
\[
\lambda\leq p-\frac{(1-B)(k-p)(A-B)(p-\alpha)^{2}}{[(1-B)(k-p)+(A-B)(p-\alpha
)]^{2}(1-\mu)^{k-p}\tfrac{\Gamma(k+\delta)}{\Gamma(p+\delta)}-[(A-B)(p-\alpha
)]^{2}}\ (k\geq p+1).
\]
\noindent Now, defining the function $\Phi(k)$ by%
\[
\Phi(k)=p-\frac{(1-B)(k-p)(A-B)(p-\alpha)^{2}}{[(1-B)(k-p)+(A-B)(p-\alpha
)]^{2}(1-\mu)^{k-p}\tfrac{\Gamma(k+\delta)}{\Gamma(p+\delta)}-[(A-B)(p-\alpha
)]^{2}}\ (k\geq p+1).
\]
\noindent We see that $\Phi(k)$\ is an increasing function of $k$\ $(k\geq
p+1)$. Therefore, we conclude that%

\begin{equation}
\lambda=\Phi(p+1)=p-\frac{(1-B)(A-B)(1-\alpha)^{2}}{[(1-B)+(A-B)(p-\alpha
)]^{2}(1-\mu)\tfrac{\Gamma(p+1+\delta)}{\Gamma(p+\delta)}-[(A-B)(p-\alpha
)]^{2}},\tag{2.27}%
\end{equation}
\noindent this completes the proof.\newline

\noindent\textbf{Theorem 6.} If $f_{1}\in R_{\mu,p}^{\delta}(\alpha;A,B)$ and
$f_{2}\in R_{\mu,p}^{\delta}(\beta;A,B)$, which are defined by (2.22), then
$(f_{1}\ast f_{2})(z)\in R_{\mu,p}^{\delta}(\xi;A,B)$, where%
\begin{equation}
\xi=p-\tfrac{(1-B)(A-B)(p-\alpha)(p-\beta)}{[(1-B)+(A-B)(p-\alpha
)][(1-B)+(A-B)(p-\beta)](1-\mu)^{k-p}\tfrac{\Gamma(k+\delta)}{\Gamma
(p+\delta)}-(A-B)^{2}(p-\alpha)(p-\beta)}\tag{2.28}%
\end{equation}

\noindent the result is the best possible for the functions%
\begin{equation}
f_{1}(z)=z^{p}-\frac{(A-B)(p-\alpha)}{[(1-B)+(A-B)(p-\alpha)](1-\mu
)\dfrac{\Gamma(p+\delta+1)}{\Gamma(p+\delta)}}z^{p+1},\tag{2.29}%
\end{equation}
\noindent and%
\begin{equation}
f_{2}(z)=z^{p}-\frac{(A-B)(p-\beta)}{[(1-B)+(A-B)(p-\beta)](1-\mu
)\dfrac{\Gamma(p+\delta+1)}{\Gamma(p+\delta)}}z^{p+1}.\tag{2.30}%
\end{equation}

\section*{\noindent3. Applications of Fractional Calculus}

In our present investigation, we shall make use of the familiar integral
operator $J_{c,p}$ defined by (see \cite{Jc 1}, \cite{Jc 2} and \cite{Jc 3})%

\begin{equation}
(J_{c,p}f)(z)=\frac{c+p}{z^{c}}%
{\displaystyle\int\limits_{0}^{z}}
t^{c-1}f(t)dt\text{ \ \ }(f\in T(p);c>-p),\tag{3.1}%
\end{equation}
\noindent also, the fractional integral of order $\eta$ is defined, for a
function $f$, by%

\begin{equation}
D_{z}^{-\eta}f(z)=\frac{1}{\Gamma(\eta)}%
{\displaystyle\int\limits_{0}^{z}}
\frac{f(\xi)}{(z-\xi)^{1-\eta}}d\xi\text{ \ \ }(\eta>0),\tag{3.2}%
\end{equation}
\noindent where the function $f$ is analytic in a simply-connected domain of
the complex plane containing the origin and the multiplicity of $(z-\xi
)^{\eta-1}$ is removed by requiring $\log(z-\xi)$ to be real when $z-\xi>0$.

\noindent The fractional derivative of order $\eta$ is defined, for a function
$f$, by%
\begin{equation}
D_{z}^{\eta}f(z)=\frac{1}{\Gamma(1-\eta)}\frac{d}{dz}%
{\displaystyle\int\limits_{0}^{z}}
\frac{f(\xi)}{(z-\xi)^{\eta}}d\xi\text{ \ \ }(0\leq\eta<1),\tag{3.3}%
\end{equation}
\noindent where the function $f(z)$ is constrained, and the multiplicity of
$(z-\xi)^{-\eta}$ is removed, as above. In this section, we investigate the
distortion properties of functions in the class $R_{\mu,p}^{\delta}%
(\alpha;A,B)$ involving the operator $J_{c,p}$\ and the fractional calculus
$D_{z}^{-\eta}$ and $D_{z}^{\eta}$. By using (3.1), (3.2), (3.3) and (1.1) it
is easily to deduce that:%
\begin{equation}
D_{z}^{\eta}\left(  J_{c,p}f(z)\right)  =\frac{\Gamma(p+1)}{\Gamma(2-\eta
)}z^{p-\eta}-%
{\textstyle\sum\limits_{k=p+1}^{\infty}}
\frac{(c+p)\Gamma(k+1)}{(c+k)\Gamma(k-\eta+1)}a_{k}z^{k-\eta},\tag{3.4}%
\end{equation}
\noindent and%
\begin{equation}
D_{z}^{-\mu}\left(  J_{c,p}f(z)\right)  =\frac{\Gamma(p+1)}{\Gamma(2+\eta
)}z^{p+\eta}-%
{\textstyle\sum\limits_{k=p+1}^{\infty}}
\frac{(c+p)\Gamma(k+1)}{(c+k)\Gamma(k+\eta+1)}a_{k}z^{k+\eta}.\tag{3.5}%
\end{equation}

\noindent\textbf{Theorem 7. }Let the function $f$ defined by (1.1) be in the
class $R_{\mu,p}^{\delta}(\alpha;A,B).$ Then%
\begin{equation}
\left\vert D_{z}^{-\eta}\left(  J_{c.p}f(z)\right)  \right\vert \geq\left\{
\tfrac{\Gamma(p+1)}{\Gamma(p+1+\eta)}-\tfrac{(c+p)\Gamma(p+2)(B-A)(p-\alpha
)}{(c+p+1)\Gamma(p+\eta+2)[(1-B)+(A-B)(p-\alpha)](1-\mu)\frac{\Gamma
(p+\delta+1)}{\Gamma(p+\delta)}}\left\vert z\right\vert \right\}  \left\vert
z\right\vert ^{p+\eta},\tag{3.6}%
\end{equation}
and%
\begin{equation}
\left\vert D_{z}^{-\eta}\left(  J_{c,p}f(z)\right)  \right\vert \leq\left\{
\tfrac{\Gamma(p+1)}{\Gamma(p+1+\eta)}+\tfrac{(c+p)\Gamma(p+2)(A-B)(p-\alpha
)}{(c+p+1)\Gamma(p-\eta+2)[(1-B)+(A-B)(p-\alpha)](1-\mu)\frac{\Gamma
(p+\delta+1)}{\Gamma(p+\delta)}}\left\vert z\right\vert \right\}  \left\vert
z\right\vert ^{p+\eta},\tag{3.7}%
\end{equation}

\noindent these results are sharp.

\noindent\textbf{Proof}. In view of Theorem 1 we have%

\[
\frac{\lbrack(1-B)+(A-B)(p-\alpha)](1-\mu)\frac{\Gamma(p+\delta+1)}%
{\Gamma(p+\delta)}}{(A-B)(p-\alpha)}%
{\textstyle\sum\limits_{k=p+1}^{\infty}}
a_{k}\leq
\]

\[%
{\textstyle\sum\limits_{k=p+1}^{\infty}}
\frac{[(1-B)+(A-B)(p-\alpha)](1-\mu)\frac{\Gamma(p+\delta+1)}{\Gamma
(p+\delta)}}{(A-B)(p-\alpha)}a_{k}\leq1,
\]
which readily yields%
\begin{equation}%
{\textstyle\sum\limits_{k=p+1}^{\infty}}
a_{k}\leq\frac{(A-B)(p-\alpha)}{[(1-B)+(A-B)(p-\alpha)](1-\mu)\frac
{\Gamma(p+\delta+1)}{\Gamma(p+\delta)}}.\tag{3.8}%
\end{equation}
Suppose that function $F(z)$ defined in $U$ by%

\begin{align}
F(z)  & =\frac{\Gamma(p+\eta+1)}{\Gamma(p+1)}z^{-\eta}\left[  D_{z}^{-\eta
}\left(  J_{c,p}f(z)\right)  \right] \nonumber\\
& =z^{p}-%
{\textstyle\sum\limits_{k=p+1}^{\infty}}
\frac{(c+p)\Gamma(k+1)\Gamma(p+\eta+1)}{(c+k)\Gamma(p+1)\Gamma(k+\eta+1)}%
a_{k}z^{k}\nonumber\\
& =z^{p}-%
{\textstyle\sum\limits_{k=p+1}^{\infty}}
\Upsilon(k)a_{k}z^{k},\tag{3.9}%
\end{align}
where
\begin{equation}
\Upsilon(k)=\frac{(c+p)\Gamma(k+1)\Gamma(p+\eta+1)}{(c+k)\Gamma(p+1)\Gamma
(k+\eta+1)}.\tag{3.10}%
\end{equation}
Since $\Upsilon(k)$ is a decreasing function of $k$,%
\begin{equation}
0<\Upsilon(k)\leq\Upsilon(p+1)=\frac{(c+p)(p+1)\Gamma(p+\eta+1)}%
{(c+p+1)(p+\eta+1)}.\tag{3.11}%
\end{equation}
By using (3.9) and (3.11) we have%
\begin{align}
\left\vert F(z)\right\vert  & =\left\vert z^{p}-%
{\textstyle\sum\limits_{k=p+1}^{\infty}}
\Upsilon(k)a_{k}z^{k}\right\vert \geq\left\vert z\right\vert ^{p}%
-\Upsilon(p+1)\left\vert z\right\vert ^{p+1}%
{\textstyle\sum\limits_{k=p+1}^{\infty}}
a_{k}\nonumber\\
& \geq\left\vert z\right\vert ^{p}-\left\{  \tfrac{(c+p)\Gamma(p+2)\Gamma
(p+\eta+1)(A-B)(p-\alpha)}{(c+p+1)\Gamma(p+1)\Gamma(p+\eta
+2)[(1-B)+(A-B)(p-\alpha)](1-\mu)\frac{\Gamma(p+\delta+1)}{\Gamma(p+\delta)}%
}\right\}  \left\vert z\right\vert ^{p+1},\tag{3.12}%
\end{align}
and%
\begin{align}
\left\vert F(z)\right\vert  & =\left\vert z^{p}+%
{\textstyle\sum\limits_{k=p+1}^{\infty}}
\Upsilon(k)a_{k}z^{k}\right\vert \leq\left\vert z\right\vert ^{p}%
-\Upsilon(p+1)\left\vert z\right\vert ^{p+1}%
{\textstyle\sum\limits_{k=p+1}^{\infty}}
a_{k}\nonumber\\
& \leq\left\vert z\right\vert ^{p}+\left\{  \tfrac{(c+p)\Gamma(p+2)\Gamma
(p+\eta+1)(A-B)(p-\alpha)}{(c+p+1)\Gamma(p+1)\Gamma(p+\eta
+2)[(1-B)+(A-B)(p-\alpha)](1-\mu)\frac{\Gamma(p+\delta+1)}{\Gamma(p+\delta)}%
}\right\}  \left\vert z\right\vert ^{p+1},\tag{3.13}%
\end{align}
which yield the inequalities (3.6) and (3.7) of Theorem 8.

\noindent This equalities in (3.6) and (3.7) are attained for the function $f
$ of which%
\begin{equation}
D_{z}^{-\eta}\left(  J_{c,p}f(z)\right)  =\left\{  \tfrac{\Gamma(p+1)}%
{\Gamma(p+1+\eta)}-\tfrac{(c+p)\Gamma(p+2)(A-B)(p-\alpha)}{(c+p+1)\Gamma
(p+1)\Gamma(p+\eta+2)[(1-B)+(A-B)(p-\alpha)](1-\mu)\frac{\Gamma(p+\delta
+1)}{\Gamma(p+\delta)}}z\right\}  z^{p+\eta},\tag{3.14}%
\end{equation}
or, equivalently%
\[
J_{c,p}f(z)=z^{p}-\tfrac{(c+p)(A-B)(p-\alpha)}{(c+p+1)[(1-B)+(A-B)(p-\alpha
)](1-\mu)\frac{\Gamma(p+\delta+1)}{\Gamma(p+\delta)}}z^{p+1},
\]
\newline Thus the proof of Theorem 7 is completed.\newline\noindent Another
inequalities can be given and the proof is omitted.\newline\textbf{Theorem 8.
}Let the function $f$ defined by (1.1) be in the class $R_{\mu,p}^{\delta
}(\alpha;A,B).$ Then%
\begin{equation}
\left\vert D_{z}^{\eta}\left(  J_{c,p}f(z)\right)  \right\vert \geq\left\{
\tfrac{\Gamma(p+1)}{\Gamma(p+1+\eta)}-\tfrac{(c+p)\Gamma(p+2)(A-B)(p-\alpha
)}{(c+p+1)\Gamma(p+1)\Gamma(p+\eta+2)[(1-B)+(A-B)(p-\alpha)](1-\mu
)\frac{\Gamma(p+\delta+1)}{\Gamma(p+\delta)}}\left\vert z\right\vert \right\}
\left\vert z\right\vert ^{p-\eta},\tag{3.15}%
\end{equation}

\noindent and%
\begin{equation}
\left\vert D_{z}^{\eta}\left(  J_{c,p}f(z)\right)  \right\vert \leq\left\{
\tfrac{\Gamma(p+1)}{\Gamma(p+1+\eta)}+\tfrac{(c+p)\Gamma(p+2)(A-B)(p-\alpha
)}{(c+p+1)\Gamma(p+1)\Gamma(p+\eta+2)[(1-B)+(A-B)(p-\alpha)](1-\mu
)\frac{\Gamma(p+\delta+1)}{\Gamma(p+\delta)}}\left\vert z\right\vert \right\}
\left\vert z\right\vert ^{p-\eta}.\tag{3.16}%
\end{equation}
\noindent Each of the assertions (3.15) and (3.16) is sharp.\newline Then, we
can easily obtain the following two theorems and the proofs are
omitted\newline\textbf{Theorem 9. }Let the function $f$ defined by (1.1) be in
the class $R_{\mu,p}^{\delta}(\alpha;A,B).$ Then%
\begin{equation}
\left\vert J_{c}\left(  D_{z}^{\eta}f(z)\right)  \right\vert \geq\left\{
\tfrac{(c+p)}{(c-\eta+1)\Gamma(p+1-\eta)}-\tfrac{(c+p)\Gamma
(p+2)(A-B)(p-\alpha)}{(c+p+1)\Gamma(p+1)\Gamma(p+\eta+2)[(1-B)+(A-B)(p-\alpha
)](1-\mu)\frac{\Gamma(p+\delta+1)}{\Gamma(p+\delta)}}\left\vert z\right\vert
\right\}  \left\vert z\right\vert ^{p-\eta},\tag{3.17}%
\end{equation}
\noindent and%
\begin{equation}
\left\vert J_{c}\left(  D_{z}^{\eta}f(z)\right)  \right\vert \leq\left\{
\tfrac{(c+p)}{(c-\eta+1)\Gamma(p+1-\eta)}-\tfrac{(c+p)\Gamma
(p+2)(A-B)(p-\alpha)}{(c+p+1)\Gamma(p+1)\Gamma(p+\eta+2)[(1-B)+(A-B)(p-\alpha
)](1-\mu)\frac{\Gamma(p+\delta+1)}{\Gamma(p+\delta)}}\left\vert z\right\vert
\right\}  \left\vert z\right\vert ^{p-\eta}.\tag{3.18}%
\end{equation}
\noindent\textbf{Theorem 10. }If the function $f$ given by (1.1) be in the
class $R_{\mu,p}^{\delta}(\alpha;A,B).$ Then%
\begin{equation}
\left\vert J_{c}\left(  D_{z}^{-\eta}f(z)\right)  \right\vert \geq\left\{
\tfrac{(c+p)}{(c+\eta+1)\Gamma(p+1+\eta)}-\tfrac{(c+p)\Gamma
(p+2)(A-B)(p-\alpha)}{(c+p+1)\Gamma(p+1)\Gamma(p+\eta+2)[(1-B)+(A-B)(p-\alpha
)](1-\mu)\frac{\Gamma(p+\delta+1)}{\Gamma(p+\delta)}}\left\vert z\right\vert
\right\}  \left\vert z\right\vert ^{p+\eta},\tag{3.19}%
\end{equation}
\noindent and%
\begin{equation}
\left\vert J_{c}\left(  D_{z}^{-\eta}f(z)\right)  \right\vert \leq\left\{
\tfrac{(c+p)}{(c+\eta+1)\Gamma(p+1+\eta)}+\tfrac{(c+p)\Gamma
(p+2)(A-B)(p-\alpha)}{(c+p+1)\Gamma(p+1)\Gamma(p+\eta+2)[(1-B)+(A-B)(p-\alpha
)](1-\mu)\frac{\Gamma(p+\delta+1)}{\Gamma(p+\delta)}}\left\vert z\right\vert
\right\}  \left\vert z\right\vert ^{p+\eta}.\tag{3.20}%
\end{equation}

\noindent\textbf{Remark }\newline\textbf{\ }By using the coefficients
estimates (given by Theorem 2) of functions belonging to the subclass
$P_{\mu,p}^{\delta}(\alpha;A,B)$ and performing the same techniques of proofs
given during Sections 2 and 3, then we can obtain the corresponding results of
$P_{\mu,p}^{\delta}(\alpha;A,B).$

\end{document}